\documentclass[11pt,twoside]{article}
\usepackage{atmp}
\usepackage{graphicx}
\usepackage{amsmath}
\usepackage{amsfonts}
\usepackage{amssymb}
\copyrightnotice{2002}{6}{575}{591}
\newtheorem{theorem}{Theorem}

\newtheorem{claim}[theorem]{Claim}

\newtheorem{definition}[theorem]{Definition}

\newtheorem{lemma}[theorem]{Lemma}

\newtheorem{proposition}[theorem]{Proposition}

\begin{document}

\title{Topological Quantum Field Theory for \\Calabi-Yau threefolds and $G_{2}$-manifolds}
\author{Naichung Conan Leung}
%\maketitle

\address{School of Mathematics, \\University of Minnesota, \\Minneapolis, MN
55454, \\USA.
}
\vspace*{.2in}

\addressemail{LEUNG@MATH.UMN.EDU}

\markboth
{TOPOLOGICAL QUANTUM FIELD THEORY ...}
{NAICHUNG CONAN LEUNG}

\section{Introduction}

In the past two decades we have witnessed many fruitful interactions between
mathematics and physics. One example is in the Donaldson-Floer theory for
oriented four manifolds. Physical considerations lead to the discovery of the
Seiberg-Witten theory which has profound impact to our understandings of four
manifolds. Another example is in the mirror symmetry for Calabi-Yau manifolds.
This duality transformation in the string theory leads to many surprising
predictions in the enumerative geometry.
\cutpage

String theory in physics studies a ten dimensional space-time $X\times
\mathbb{R}^{3,1}$. Here $X$ a six dimensional Riemannian manifold with its
holonomy group inside $SU\left(  3\right)  $, the so-called \textit{Calabi-Yau
threefold}. Certain parts of the mirror symmetry conjecture, as studied by
Vafa's group, are specific for Calabi-Yau manifolds of complex dimension
\textit{three}. They include the Gopakumar-Vafa conjecture for the
Gromov-Witten invariants of \textit{arbitrary} genus, the Ooguri-Vafa
conjecture on the relationships between knot invariants and enumerations of
holomorphic disks and so on. The key reason is they belong to a duality
theory for $G_{2}$-manifolds. $G_{2}$-manifolds can be naturally interpreted
as special Octonion manifolds \cite{Le RG over A}. For any Calabi-Yau
threefold $X$, the seven dimensional manifold $X\times S^{1}$ is automatically
a $G_{2}$-manifold because of the natural inclusion $SU\left(  3\right)
\subset G_{2}$.

In recent years, there are many studies of $G_{2}$-manifolds in M-theory
including works of Archaya, Atiyah, Gukov, Vafa, Witten, Yau, Zaslow and many
others (e.g. \cite{Acharya}, \cite{Atiyah Witten}, \cite{GYZ}, \cite{Mina
Vafa}).

In the studies of the symplectic geometry of a Calabi-Yau threefold $X$, we
consider unitary flat bundles over three dimensional (special) Lagrangian
submanifolds $L$ in $X$. The corresponding geometry for a $G_{2}$-manifold $M
$ is called the \textit{special }$\mathbb{H}$\textit{-Lagrangian geometry} (or
\textit{C-geometry }in \cite{LL}). where we consider Anti-Self-Dual (abbrev.
ASD) bundles over four dimensional coassociative submanifolds, or equivalently
\textit{special }$\mathbb{H}$\textit{-Lagrangian submanifolds of type II
}\cite{Le RG over A},\textit{\ } (abbrev. $\mathbb{H}$-SLag) $C$ in $M$.

Counting ASD bundles over a fixed four manifold $C$ is the well-known theory
of Donaldson differentiable invariants, $Don\left(  C\right)  $. Similarly,
counting unitary flat bundles over a fixed three manifold $L$ is Floer's
Chern-Simons homology theory, $HF_{CS}\left(  L\right)  $. When $C$ is a
connected sum $C_{1}\#_{L}C_{2}$ along a homology three sphere, the relative
Donaldson invariants $Don\left(  C_{i}\right)  $'s take values in
$HF_{CS}\left(  L\right)  $ and $Don\left(  C\right)  $ can be recovered from
individual pieces by a gluing theorem, $Don\left(  C\right)  =\left\langle
Don\left(  C_{1}\right)  ,Don\left(  C_{2}\right)  \right\rangle
_{HF_{CS}\left(  L\right)  }$ (see e.g. \cite{Don Instanton Book}). Similarly
when $L$ has a handlebody decomposition $L=L_{1}\#_{\Sigma}L_{2}$, each
$L_{i}$ determines a Lagrangian subspace $\mathcal{L}_{i}$ in the moduli space
$\mathcal{M}^{flat}\left(  \Sigma\right)  $ of unitary flat bundles over the
Riemann surface $\Sigma$ and Atiyah conjectures that we can recover
$HF_{CS}\left(  L\right)  $ from the Floer's Lagrangian intersection homology
group of $\mathcal{L}_{1}$ and $\mathcal{L}_{2}$ in $\mathcal{M}^{flat}\left(
\Sigma\right)  $, $HF_{CS}\left(  L\right)  =HF_{Lag}^{\mathcal{M}%
^{flat}\left(  \Sigma\right)  }\left(  \mathcal{L}_{1},\mathcal{L}_{2}\right)
$. Such algebraic structures in the Donaldson-Floer theory can be formulated
as a Topological Quantum Field Theory (abbrev. TQFT), as defined by Segal and
Atiyah \cite{At 3 4}.

In this paper, we propose a construction of a TQFT by counting ASD bundles
over four dimensional $\mathbb{H}$-SLag $C$ in any closed (almost) $G_{2}%
$-manifold $M$. We call these $\mathbb{H}$\textit{-SLag cycles} and they can
be identified as zeros of a naturally defined closed one form on the
configuration space of topological cycles. We expect to obtain a homology
theory $H_{C}\left(  M\right)  $ by applying the construction in the Witten's
Morse theory. When $M$ is non-compact with an asymptotically cylindrical end,
$X\times\lbrack0,\infty)$, then the collection of boundary data of relative
$\mathbb{H}$-SLag cycles determines a Lagrangian submanifold $\mathcal{L}_{M}$
in the moduli space $\mathcal{M}^{SLag}\left(  X\right)  $ of special
Lagrangian cycles in the Calabi-Yau threefold $X$.

When we decompose $M=M_{1}\#_{X}M_{2}$ along an infinite asymptotically
cylindrical neck, it is reasonable to expect to have a gluing formula,
\[
H_{C}\left(  M\right)  =HF_{Lag}^{\mathcal{M}^{SLag}\left(  X\right)  }\left(
\mathcal{L}_{M_{1}},\mathcal{L}_{M_{2}}\right)  \text{.}%
\]
The main technical difficulty in defining this TQFT rigorously is the
\textit{compactness} issue for the moduli space of $\mathbb{H}$-SLag cycles in
$M$. We do not know how to resolve this problem and our homology groups are
only defined in the \textit{formal} sense (and physical sense?).

%\cutpage
\section{$G_{2}$-manifolds and $\mathbb{H}$-SLag geometry}

We first review some basic definitions and properties of $G_{2}$-geometry, see
\cite{LL} for more details.

\begin{definition}
A seven dimensional Riemannian manifold $M$ is called a $G_{2}$-manifold if
the holonomy group of its Levi-Civita connection is inside $G_{2}\subset
SO\left(  7\right)  $.
\end{definition}

The simple Lie group $G_{2}$ can be identified as the subgroup of $SO\left(
7\right)  $ consisting of isomorphism $g:\mathbb{R}^{7}\rightarrow
\mathbb{R}^{7}$ preserving the linear three form $\Omega$,
\[
\Omega=f^{1}f^{2}f^{3}-f^{1}\left(  e^{1}e^{0}+e^{2}e^{3}\right)
-f^{2}\left(  e^{2}e^{0}+e^{3}e^{1}\right)  -f^{3}\left(  e^{3}e^{0}%
+e^{1}e^{2}\right)  \text{,}%
\]
where $e^{0},e^{1},e^{2},e^{3},f^{1},f^{2},f^{3}$ is any given orthonormal
frame of $\mathbb{R}^{7}$. Such a three form, or up to conjugation by elements
in $GL\left(  7,\mathbb{R}\right)  $, is called \textit{positive}, and it
determines a unique compatible inner product on $\mathbb{R}^{7}$ \cite{Bryant
Metric Exc}.

Gray \cite{Gray VectorCrossProd} shows that $G_{2}$-holonomy of $M$ can be
characterized by the existence of a positive harmonic three form $\Omega$.

\begin{definition}
A seven dimensional manifold $M$ equipped with a positive closed three form
$\Omega$ is called an almost $G_{2}$-manifold.
\end{definition}

Remark: The relationship between $G_{2}$-manifolds and almost $G_{2}%
$-manifolds is analogous to the relationship between Kahler manifolds and
symplectic manifolds. Namely we replace a parallel non-degenerate form by a
closed one.

For example, suppose that $X$ is a complex three dimensional K\"{a}hler
manifold with a trivial canonical line bundle, i.e. there exists a
nonvanishing holomorphic three form $\Omega_{X}$. Yau's celebrated theorem
says that there is a K\"{a}hler form $\omega_{X}$ on $X$ such that the
corresponding Kahler metric has holonomy in $SU\left(  3\right)  $, i.e. a
Calabi-Yau threefold. In particular both $\Omega_{X}$ and $\omega_{X}$ are
parallel forms. Then the product $M=X\times S^{1}$ is a $G_{2}$-manifold with
\[
\Omega=\operatorname{Re}\Omega_{X}+\omega_{X}\wedge d\theta\text{.}%
\]
Conversely, one can prove, using Bochner arguments, every $G_{2}$-metric on
$X\times S^{1}$ must be of this form. More generally, if $\omega_{X}$ is a
general K\"{a}hler form on $X$, then $\left(  X\times S^{1},\Omega\right)  $
is an \textit{almost} $G_{2}$-manifold and the converse is also true.

\bigskip

Next we quickly review the geometry of $\mathbb{H}$-SLag cycles in an almost
$G_{2}$-manifold (see \cite{LL}).

\begin{definition}
An orientable four dimensional submanifold $C$ in an almost $G_{2}$-manifold
$\left(  M,\Omega\right)  $ is called a coassociative submanifold, or simply a
$\mathbb{H}$-SLag, if the restriction of $\Omega$ to $C$ is identically zero,
\[
\Omega|_{C}=0\text{.}%
\]
\end{definition}

If $M$ is a $G_{2}$-manifold, then any coassociative submanifold $C$ in $M$ is
calibrated by $\ast\Omega$ in the sense of Harvey and Lawson \cite{HL}, in
particular, it is an absolute minimal submanifold in $M$. The normal bundle of
any $\mathbb{H}$-SLag $C$ can be naturally identified with the bundle of
self-dual two forms on $C$. McLean \cite{McLean} shows that infinitesimal
deformations of any $\mathbb{H}$-SLag are unobstructed and they are
parametrized by the space of harmonic self-dual two forms on $C$, i.e.
$H_{+}^{2}\left(  C,\mathbb{R}\right)  $.

For example, if $S$ is a complex surface in a Calabi-Yau threefold $X$, then
$S\times\left\{  t\right\}  $ is a $\mathbb{H}$-SLag in $M=X\times S^{1}$ for
any $t\in S^{1}$. Notice that $H_{+}^{2}\left(  S,\mathbb{R}\right)  $ is
spanned by the Kahler form and the real and imaginary parts of holomorphic two
forms on $S$, and the latter can be identified holomorphic normal vector
fields along $S$ because of the adjunction formula and the Calabi-Yau
condition on $X$. Thus all deformations of $S\times\left\{  t\right\}  $ in $M
$ as $\mathbb{H}$-SLag submanifolds are of the same form. Similarly, if $L$ is
a three dimensional special Lagrangian submanifold in $X$ with phase $\pi/2$,
i.e. $\omega|_{L}=\operatorname{Re}\Omega_{X}|_{L}=0$, then $L\times S^{1}$ is
also a $\mathbb{H}$-SLag in $M=X\times S^{1}$. Furthermore, all deformations
of $L\times S^{1}$ in $M$ as $\mathbb{H}$-SLag submanifolds are of the same
form because $H_{+}^{2}\left(  L\times S^{1}\right)  \cong H^{1}\left(
L\right)  $, which parametrizes infinitesimal deformations of special
Lagrangian submanifolds in $X$.

\begin{definition}
A $\mathbb{H}$-SLag cycle in an almost $G_{2}$-manifold $\left(
M,\Omega\right)  $ is a pair $\left(  C,D_{E}\right)  $ with $C$ a
$\mathbb{H}$-SLag in $M$ and $D_{E}$ an ASD connection over $C$.
\end{definition}

Remark: $\mathbb{H}$-SLag cycles are supersymmetric cycles in physics as
studied in \cite{MMMS}. Their moduli space admits a natural three form and a
cubic tensor \cite{LL}, which play the roles of the correlation function and
the Yukawa coupling in physics.

We assume that the ASD connection $D_{E}$ over $C$ has rank one, i.e. a
$U\left(  1\right)  $ connection. This avoids the occurrence of reducible
connections, thus the moduli space $\mathcal{M}^{\mathbb{H}-SLag}\left(
M\right)  $ of $\mathbb{H}$-SLag cycles in $M$ is a smooth manifold. It has a
natural orientation and its expected dimension equals $b^{1}\left(  C\right)
$, the first Betti number of $C$. This is because the moduli space of
$\mathbb{H}$-SLags has dimension equals $b_{+}^{2}\left(  C\right)  $
\cite{McLean} and the existence of an ASD $U\left(  1\right)  $-connection
over $C$ is equivalent to $H_{-}^{2}\left(  C,\mathbb{R}\right)  \cap
H^{2}\left(  C,\mathbb{Z}\right)  \neq\phi$. The number $b^{1}\left(
C\right)  $ is responsible for twisting by a flat $U\left(  1\right)  $-connection.

\bigskip

For simplicity, we assume that $b^{1}\left(  C\right)  =0$, otherwise, one can
cut down the dimension of $\mathcal{M}^{\mathbb{H}-SLag}\left(  M\right)  $ to
zero by requiring the ASD connections over $C$ to have trivial holonomy around
loops $\gamma_{1},...,\gamma_{b^{1}\left(  C\right)  }$ in $C$ representing an
integral basis of $H_{1}\left(  C,\mathbb{Z}\right)  $. We plan to count the
algebraic number of points in this moduli space $\#\mathcal{M}^{\mathbb{H}%
-SLag}\left(  M\right)  $.

This number, in the case of $X\times S^{1}$, can be identified with a proposed
invariant of Joyce \cite{Joyce Count SLag} defined by counting rigid special
Lagrangian submanifolds in any Calabi-Yau threefold. To explain this, we need
the following proposition on the strong rigidity of product $\mathbb{H}$-SLags.

\begin{proposition}
If $L\times S^{1}$ is a $\mathbb{H}$-SLag in $M=X\times S^{1}$ with $X$ a
Calabi-Yau threefold, then any $\mathbb{H}$-SLag representing the same
homology class must also be a product.
\end{proposition}

Proof: For simplicity we assume that the volume of the $S^{1}$ factor is
unity, $Vol\left(  S^{1}\right)  =1$. If $L\times S^{1}$ is a $\mathbb{H}%
$-SLag in $M$ then $L$ is special Lagrangian submanifold in $X$ with phase
$\pi/2$, i.e. $\operatorname{Re}\Omega_{X}|_{L}=\omega|_{L}=0$. Suppose $C$ is
another $\mathbb{H}$-SLag in $M$ representing the same homology class, we have
$Vol\left(  C\right)  =Vol\left(  L\right)  $. If we write $C_{\theta}%
=C\cap\left(  X\times\left\{  \theta\right\}  \right)  $ for any $\theta\in
S^{1}$, then $Vol\left(  C_{\theta}\right)  \geq Vol\left(  L\right)  , $ as
$L$ is a calibrated submanifold in $X$. Furthermore the equality sign holds
only if $C_{\theta}$ is also calibrated. In general we have
\[
Vol\left(  C\right)  \geq\int_{S^{1}}Vol\left(  C_{\theta}\right)  d\theta,
\]
with the equality sign holds if and only if $C$ is a product with $S^{1}$.
Combining these, we have
\[
Vol\left(  L\right)  =Vol\left(  C\right)  \geq\int_{S^{1}}Vol\left(
C_{\theta}\right)  d\theta\geq\int_{S^{1}}Vol\left(  L\right)  d\theta
=Vol\left(  L\right)  \text{.}%
\]
Thus both inequalities are indeed equal. Hence $C=L^{\prime}\times S^{1}$ for
some special Lagrangian submanifold $L^{\prime}$ in $X$. $\blacksquare$

\bigskip\ 

Suppose $M=X\times S^{1}$ is a product $G_{2}$-manifold and we consider
product $\mathbb{H}$-SLag $C=L\times S^{1}$ in $M$. From the above
proposition, every $\mathbb{H}$-SLag representing $\left[  C\right]  $ must
also be a product. Since $b_{+}^{2}\left(  C\right)  =b^{1}\left(  L\right)
$, the rigidity of the $\mathbb{H}$-SLag $C$ in $M$ is equivalent to the
rigidity of the special Lagrangian submanifold $L$ in $X$. When this happens,
i.e. $L$ is a rational homology three sphere, we have $b^{2}\left(  C\right)
=0$ and
\[
\text{No. of ASD U(1)-bdl/}C=\#H^{2}\left(  C,\mathbb{Z}\right)
=\#H^{2}\left(  L,\mathbb{Z}\right)  =\#H_{1}\left(  L,\mathbb{Z}\right)
\text{.}%
\]
Here we have used the fact that the first cohomology group is always torsion
free. Thus the number of such $\mathbb{H}$-SLag cycles in $X\times S^{1}$
equals the number of special Lagrangian rational homology three spheres in a
Calabi-Yau threefold $X$, weighted by $\#H_{1}\left(  L,\mathbb{Z}\right)  $.
Joyce \cite{Joyce Count SLag} shows that with this particular weight, the
numbers of special Lagrangians in any Calabi-Yau threefold behave well under
various surgeries on $X$, and expects them to be invariants. Thus in this
case, we have
\[
\#\mathcal{M}^{\mathbb{H}-SLag}\left(  X\times S^{1}\right)  =\text{Joyce's
proposed invariant for }\#\text{SLag. in }X\text{.}%
\]
In the next section, we will propose a homology theory, whose Euler
characteristic gives $\#\mathcal{M}^{\mathbb{H}-SLag}\left(  M\right)  $.

\section{Witten's Morse theory for $\mathbb{H}$-SLag cycles}

We are going to use the parametrized version of $\mathbb{H}$-SLag cycles in
any almost $G_{2}$-manifold $M$. We fix an oriented smooth four dimensional
manifold $C$ and a rank $r$ Hermitian vector bundle $E$ over $C$. We consider
the \textit{configuration space}
\[
\mathcal{C}=Map\left(  C,M\right)  \times\mathcal{A}\left(  E\right)  ,
\]
where $\mathcal{A}\left(  E\right)  $ is the space of Hermitian connections on
$E$.

\begin{definition}
An element $\left(  f,D_{E}\right)  $ in $\mathcal{C}$ is called a
\textit{parametrized $\mathbb{H}$-SLag cycles} in $M$ if
\[
f^{\ast}\Omega=F_{E}^{+}=0\text{,}%
\]
where the self-duality is defined using the pullback metric from $M$.
\end{definition}

Instead of $Aut\left(  E\right)  $, the symmetry group $\mathcal{G}$ in our
situation consists of gauge transformations of $E$ which cover
\textit{arbitrary} diffeomorphisms on $M$,
\[%
\begin{array}
[c]{ccc}%
E & \overset{g}{\rightarrow} & E\\
\downarrow &  & \downarrow\\
M & \overset{g_{M}}{\rightarrow} & M.
\end{array}
\]
It fits into the following exact sequence,%

\[
1\rightarrow Aut\left(  E\right)  \rightarrow\mathcal{G}\rightarrow
Diff\left(  C\right)  \rightarrow1\text{.}%
\]
The natural action of $\mathcal{G}$ on $\mathcal{C}$ is given by
\[
g\cdot\left(  f,D_{E}\right)  =\left(  f\circ g_{M},g^{\ast}D_{E}\right)  ,
\]
for any $\left(  f,D_{E}\right)  \in\mathcal{C}=Map\left(  C,M\right)
\times\mathcal{A}\left(  E\right)  $. Notice that $\mathcal{G}$ preserves the
set of parametrized $\mathbb{H}$-SLag cycles in $M$.

The configuration space $\mathcal{C}$ has a natural one form $\Phi_{0}$: At
any $\left(  f,D_{E}\right)  \in\mathcal{C}$ we can identify the tangent space
of $\mathcal{C}$ as
\[
T_{\left(  f,D_{E}\right)  }\mathcal{C}=\Gamma\left(  C,f^{\ast}T_{M}\right)
\times\Omega^{1}\left(  C,ad\left(  E\right)  \right)  \text{.}%
\]
We define
\[
\Phi_{0}\left(  f,D_{E}\right)  \left(  v,B\right)  =\int_{C}Tr\left[
f^{\ast}\left(  \iota_{v}\Omega\right)  \wedge F_{E}+f^{\ast}\Omega\wedge
B\right]  \text{,}%
\]
for any $\left(  v,B\right)  \in T_{\left(  f,D_{E}\right)  }\mathcal{C}$.

\begin{proposition}
The one form $\Phi_{0}$ on $\mathcal{C}$ is closed and invariant under the
action by $\mathcal{G}$.
\end{proposition}

Proof: Recall that there is a universal connection $\mathbb{D}_{E}$ over
$C\times\mathcal{A}\left(  E\right)  $ whose curvature $\mathbb{F}_{E}$ at a
point $\left(  x,D_{E}\right)  $ equals,
\begin{align*}
\mathbb{F}_{E}|_{\left(  x,D_{E}\right)  }  &  =\left(  \mathbb{F}_{E}%
^{2,0},\mathbb{F}_{E}^{1,1},\mathbb{F}_{E}^{0,2}\right) \\
&  \in\Omega^{2}\left(  C\right)  \otimes\Omega^{0}\left(  \mathcal{A}\right)
+\Omega^{1}\left(  C\right)  \otimes\Omega^{1}\left(  \mathcal{A}\right)
+\Omega^{0}\left(  C\right)  \otimes\Omega^{2}\left(  \mathcal{A}\right)
\end{align*}
with
\[
\mathbb{F}_{E}^{2,0}=F_{E},\,\mathbb{F}_{E}^{1,1}\left(  v,B\right)  =B\left(
v\right)  ,\,\mathbb{F}_{E}^{0,2}=0,
\]
where $v\in T_{x}C$ and $B\in\Omega^{1}\left(  C,ad\left(  E\right)  \right)
=T_{D_{E}}\mathcal{A}\left(  E\right)  $ (see e.g. \cite{Le Sympl Gauge}). The
Bianchi identity implies that $Tr\mathbb{F}_{E}$ is a closed form on
$C\times\mathcal{A}\left(  E\right)  $. We also consider the evaluation map,
\begin{gather*}
ev:C\times Map\left(  C,M\right)  \rightarrow M\\
ev\left(  x,f\right)  =f\left(  x\right)  .
\end{gather*}
It is not difficult to see that the pushforward of the differential form
$ev^{\ast}\left(  \Omega\right)  \wedge Tr\mathbb{F}_{E}$ on $C\times
Map\left(  C,M\right)  \times\mathcal{A}\left(  E\right)  $ to $Map\left(
C,M\right)  \times\mathcal{A}\left(  E\right)  $ equals $\Phi_{0}$, i.e.
\[
\Phi_{0}=\int_{C}ev^{\ast}\left(  \Omega\right)  \wedge Tr\mathbb{F}%
_{E}\text{.}%
\]
Therefore the closedness of $\Phi_{0}$ follows from the closedness of $\Omega
$. It is also clear from this description of $\Phi_{0}$ that it is
$\mathcal{G}$-invariant. $\blacksquare$

\bigskip

From this proposition, we know that $\Phi_{0}=d\Psi_{0}$ locally for some
function $\Psi_{0}$ on $\mathcal{C}$. As in the Chern-Simons theory, this
function $\Psi_{0}$ can be obtained explicitly by integrating the closed one
form $\Phi_{0}$ along any path joining to a fixed element in $\mathcal{C}$.
When $M=X\times S^{1}$ and $C=L\times S^{1}$, this is essentially the
functional used by Thomas in \cite{Thomas}.

From now on, we assume that $E$ is a rank one bundle.

\begin{lemma}
The zeros of $\Phi_{0}$ are the same as parametrized $\mathbb{H}$-SLag cycles
in $M$.
\end{lemma}

Proof: Suppose $\left(  f,D_{E}\right)  $ is a zero of $\Phi_{0}$. By
evaluating it on various $\left(  0,B\right)  $, we have $f^{\ast}\Omega=0$,
i.e. $f:C\rightarrow M$ is a parametrized $\mathbb{H}$-SLag. This implies that
the map
\[
\lrcorner\Omega:T_{f\left(  x\right)  }M\rightarrow\Lambda^{2}T_{x}^{\ast}C
\]
has image equals $\Lambda_{+}^{2}T_{x}^{\ast}C$, for any $x\in C$. By
evaluating $\Phi_{0}$ on various $\left(  v,0\right)  $, we have $F_{E}^{+}%
=0$, i.e. $\left(  f,D_{E}\right)  $ is a parametrized $\mathbb{H}$-SLag cycle
in $M$. The converse is obvious. $\blacksquare$

\bigskip

From above results, $\Phi_{0}$ descends to a closed one form on $\mathcal{C}%
/\mathcal{G}$, called $\Phi$. Locally we can write $\Phi=d\mathcal{F}$ for
some function $\mathcal{F}$ whose critical points are precisely
(unparametrized) $\mathbb{H}$-SLag cycles in $M$. Using the gradient flow
lines of $\mathcal{F}$, we could formally define a Witten's Morse homology
group, as in the famous Floer's theory. Roughly speaking one defines a complex
$\left(  \mathbf{C}_{\ast},\partial\right)  $, where $\mathbf{C}_{\ast} $ is
the free Abelian group generated by critical points of $\mathcal{F}$ and
$\partial$ is defined by counting the number of gradient flow lines between
two critical points of relative index one.

Remark: The equations for the gradient flow are given by
\[
\frac{\partial f}{\partial t}=\ast\left(  f^{\ast}\xi\wedge F_{E}\right)
,\,\frac{\partial D_{E}}{\partial t}=\ast\left(  f^{\ast}\Omega\right)  ,
\]
where $\xi\in\Omega^{2}\left(  M,T_{M}\right)  $ is defined by $\left\langle
\xi\left(  u,v\right)  ,w\right\rangle =\Omega\left(  u,v,w\right)  $.

The equation
\[
\partial^{2}=0
\]
requires a good compactification of the moduli space of $\mathbb{H}$-SLag
cycles in $M$, which we are lacking at this moment (see \cite{Tian} however).
We denote this proposed homology group as $H_{C}\left(  M\right)  $, or
$H_{C}\left(  M,\alpha\right)  $ when $f_{\ast}\left[  C\right]  =\alpha\in
H_{4}\left(  M,\mathbb{Z}\right)  $.

This homology group should be invariant under deformations of the almost
$G_{2}$-metric on $M$ and its Euler characteristic equals,
\[
\chi\left(  H_{C}\left(  M\right)  \right)  =\#\mathcal{M}^{\mathbb{H}%
-SLag}\left(  M\right)  \text{.}%
\]
Like Floer homology groups, they measure the \textit{middle dimensional}
topology of the configuration space $\mathcal{C}$ divided by $\mathcal{G}$.

\section{TQFT of $\mathbb{H}$-SLag cycles}

In this section we study complete almost $G_{2}$-manifold $M_{i}$ with
asymptotically cylindrical ends and the behavior of $H_{C}\left(  M\right)  $
when a closed almost $G_{2}$-manifold $M$ decomposes into connected sum of two
pieces, each with an asymptotically cylindrical end,
\[
M=M_{1}\underset{X}{\#}M_{2}.
\]
Nontrivial examples of compact $G_{2}$-manifolds are constructed by Kovalev
\cite{Kovalev} using such connected sum approach. The boundary manifold $X$ is
necessary a Calabi-Yau threefold. We plan to discuss analytic aspects of
$M_{i}$'s in a future paper \cite{Le Asym Cyl G2}.

Each $M_{i}$'s will define a Lagrangian subspace $\mathcal{L}_{M_{i}}$ in the
moduli space of special Lagrangian cycles in $X$. Furthermore we expect to
have a gluing formula expressing the above homology group for $M$ in terms of
the Floer Lagrangian intersection homology group for the two Lagrangian
subspaces $\mathcal{L}_{M_{1}}$ and $\mathcal{L}_{M_{2}}$,
\[
H_{C}\left(  M\right)  =HF_{Lag}^{\mathcal{M}^{SLag}\left(  X\right)  }\left(
\mathcal{L}_{M_{1}},\mathcal{L}_{M_{2}}\right)  \text{.}%
\]
These properties can be reformulated to give us a topological quantum field
theory. To begin we have the following definition.

\begin{definition}
An almost $G_{2}$-manifold $M$ is called cylindrical if $M=X\times
\mathbb{R}^{1}$ and its positive three form respect such product structure,
i.e.
\[
\Omega_{0}=\operatorname{Re}\Omega_{X}+\omega_{X}\wedge dt\text{.}%
\]

A complete almost $G_{2}$-manifold $M$ with one end $X\times\lbrack0,\infty)$
is called asymptotically cylindrical if the restriction of its positive three
form equals to the above one for large $t$, up to a possible error of order
$O\left(  e^{-t}\right)  $. More precisely the positive three form $\Omega$ of
$M$ restricted to its end equals,
\[
\Omega=\Omega_{0}+d\zeta
\]
for some two form $\zeta$ satisfying $\left|  \zeta\right|  +\left|
\nabla\zeta\right|  +\left|  \nabla^{2}\zeta\right|  +\left|  \nabla^{3}%
\zeta\right|  \leq Ce^{-t}.$
\end{definition}

Remark: If $M$ is an almost $G_{2}$-manifold with an asymptotically
cylindrical end $X\times\lbrack0,\infty)$, then $\left(  X,\omega_{X}%
,\Omega_{X}\right)  $ is a complex threefold with a trivial canonical line
bundle, but the K\"{a}hler form $\omega_{X}$ might not be Einstein. This is
so, i.e. a Calabi-Yau threefold, provided that $M$ is a $G_{2}$-manifold. We
will simply write $\partial M=X$.

We consider $\mathbb{H}$-SLags $C$ in $M$ which satisfy a \textit{Neumann
condition} at infinity. That is, away from some compact set in $M$, the
immersion $f:C\rightarrow M$ can be written as
\[
f:L\times\lbrack0,\infty)\rightarrow X\times\lbrack0,\infty)
\]
with $\partial f/\partial t$ vanishes at infinite \cite{Le Asym Cyl G2}. A
relative $\mathbb{H}$-SLag itself has asymptotically cylindrical end
$L\times\lbrack0,\infty)$ with $L$ a special Lagrangian submanifold in $X$. A
\textit{relative $\mathbb{H}$-SLag cycle} in $M$ is a pair $\left(
C,D_{E}\right)  $ with $C$ a relative $\mathbb{H}$-SLag in $M$ and $D_{E}$ a
unitary connection over $C$ with finite energy,
\[
\int_{C}\left|  F_{E}\right|  ^{2}dv<\infty\text{.}%
\]
Any finite energy connection $D_{E}$ on $C$ induces a unitary flat connection
$D_{E^{\prime}}$ on $L$ \cite{Don Instanton Book}.

\bigskip

Such a pair $\left(  L,D_{E^{\prime}}\right)  $ of a unitary flat connection
$D_{E^{\prime}}$ over a special Lagrangian submanifold $L$ in a Calabi-Yau
threefold $X$ is called a \textit{special Lagrangian cycle }in $X$. Their
moduli space $\mathcal{M}^{SLag}\left(  X\right)  $ plays an important role in
the Strominger-Yau-Zaslow Mirror Conjecture \cite{SYZ} or \cite{Le Geom MS}.
The tangent space to $\mathcal{M}^{SLag}\left(  X\right)  $ is naturally
identified with $H^{2}\left(  L,\mathbb{R}\right)  \times H^{1}\left(
L,ad\left(  E^{\prime}\right)  \right)  $. For line bundles over $L$, the cup
product
\[
\cup:H^{2}\left(  L,\mathbb{R}\right)  \times H^{1}\left(  L,\mathbb{R}%
\right)  \rightarrow\mathbb{R},
\]
induces a symplectic structure on $\mathcal{M}^{SLag}\left(  X\right)  $
\cite{Hitchin SLag}. Using analytic results from \cite{Le Asym Cyl G2} about
asymptotically cylindrical manifolds, we can prove the following theorem.

\begin{claim}
Suppose $M$ is an asymptotically cylindrical (almost) $G_{2}$-manifold with
$\partial M=X$. Let $\mathcal{M}^{\mathbb{H}-SLag}\left(  M\right)  $ be the
moduli space of rank one relative $\mathbb{H}$-SLag cycles in $M$. Then the
map defined by the boundary values,
\[
b:\mathcal{M}^{\mathbb{H}-SLag}\left(  M\right)  \rightarrow\mathcal{M}%
^{SLag}\left(  X\right)  ,
\]
is a Lagrangian immersion.
\end{claim}

Sketch of the proof (\cite{Le Asym Cyl G2}): For any closed Calabi-Yau
threefold $X$ (resp. $G_{2}$-manifold $M$), the moduli space of rank one
special Lagrangian submanifolds $L$ (resp. $\mathbb{H}$-SLags $C$) is smooth
\cite{McLean} and has dimension $b^{2}\left(  L\right)  $ (resp. $b_{+}%
^{2}\left(  C\right)  $). The same holds true for complete manifold $M$ with a
asymptotically cylindrical end $X\times\lbrack0,\infty)$, where $b_{+}%
^{2}\left(  C\right)  _{L^{2}}$ denote the dimension of $L^{2}$-harmonic
self-dual two forms on a relative $\mathbb{H}$-SLag $C$ in $M$.

The linearization of the boundary value map $\mathcal{M}^{\mathbb{H}%
-SLag}\left(  M\right)  \rightarrow\mathcal{M}^{SLag}\left(  X\right)  $ is
given by $H_{+}^{2}\left(  C\right)  _{L^{2}}\overset{\alpha}{\rightarrow
}H^{2}\left(  L\right)  $. Similar for the connection part, where the boundary
value map is given by $H^{1}\left(  C\right)  _{L^{2}}\overset{\beta
}{\rightarrow}H^{1}\left(  L\right)  $. We consider the following diagram
where each row is a long exact sequence of $L^{2}$-cohomology groups for the
pair $\left(  C,L\right)  $ and each column in a perfect pairing.%

\[%
\begin{array}
[c]{cccccccccccc}%
0 & \rightarrow &  H_{+}^{2}\left(  C,L\right)  & \rightarrow &  H_{+}%
^{2}\left(  C\right)  & \overset{\alpha}{\rightarrow} & H^{2}\left(  L\right)
& \rightarrow &  H^{3}\left(  C,L\right)  & \rightarrow & \cdots & \\
&  & \otimes &  & \otimes &  & \otimes &  & \otimes &  &  & \\
0 & \leftarrow &  H_{+}^{2}\left(  C\right)  & \leftarrow &  H_{+}^{2}\left(
C,L\right)  & \leftarrow &  H^{1}\left(  L\right)  & \overset{\beta
}{\leftarrow} & H^{1}\left(  C\right)  & \leftarrow & \cdots & \\
&  & \downarrow &  & \downarrow &  & \downarrow &  & \downarrow &  &  & \\
&  & \mathbb{R} &  & \mathbb{R} &  & \mathbb{R} &  & \mathbb{R} &  &  &
\end{array}
\]
Notice that $H_{+}^{2}\left(  C,L\right)  $, $H_{+}^{2}\left(  C\right)  $ and
$H^{2}\left(  L\right)  $ parametrize infinitesimal deformation of $C$ with
fixed $\partial C$, deformation of $C$ alone and deformation of $L$ respectively.

By simple homological algebra, it is not difficult to see that
$\operatorname{Im}\alpha\oplus\operatorname{Im}\beta$ is a Lagrangian subspace
of $H^{2}\left(  L\right)  \oplus H^{1}\left(  L\right)  $ with the canonical
symplectic structure. Hence the result. $\blacksquare$

\bigskip

Remark: The deformation theory of \textit{conical} special Lagrangian
submanifolds is developed by Pacini in \cite{Pacini}.

We denote the immersed Lagrangian submanifold $b\left(  \mathcal{M}%
^{\mathbb{H}-SLag}\left(  M\right)  \right)  $ in $\mathcal{M}^{SLag}\left(
X\right)  $ by $\mathcal{L}_{M}$. When $M$ decompose as a connected sum
$M_{1}\#_{X}M_{2}$ along a long neck, as in Atiyah's conjecture on Floer
Chern-Simons homology group \cite{At 3 4}, we expect to have an isomorphism,
\[
H_{C}\left(  M\right)  \cong HF_{Lag}^{\mathcal{M}^{SLag}\left(  X\right)
}\left(  \mathcal{L}_{M_{1}},\mathcal{L}_{M_{2}}\right)  \text{.}%
\]
More precisely, suppose $\Omega_{t}$ with $t\in\lbrack0,\infty)$, is a family
of $G_{2}$-structure on $M_{t}=M$ such that as $t$ goes to infinite, $M$
decomposes into two components $M_{1}$ and $M_{2}$, each has an aymptotically
cylindrical end $X\times\lbrack0,\infty)$. Then we expect that $\lim
_{t\rightarrow\infty}H_{C}\left(  M_{t}\right)  \cong HF_{Lag}^{\mathcal{M}%
^{SLag}\left(  X\right)  }\left(  \mathcal{L}_{M_{1}},\mathcal{L}_{M_{2}%
}\right)  $. We summarize these structures in the following table:%

\[%
\begin{tabular}
[c]{|c||c|c|c}\cline{1-3}%
{\footnotesize 
Manifold:} & 
{\footnotesize 
(almost)} $G_{2}$
{\footnotesize 
-manifold,} $M^{7}$ & 
{\footnotesize 
(almost) CY threefold,
}$X^{6}$ & $%
\begin{array}
[c]{l}%
\,\\
\,
\end{array}
$\\\cline{1-3}%
{\footnotesize 
SUSY Cycles:} & $\mathbb{H}$
{\footnotesize 
-SLag. submfds.+ ASD bdl} & 
{\footnotesize 
SLag submfds.+ flat
bdl} & $%
\begin{array}
[c]{l}%
\,\\
\,
\end{array}
$\\\cline{1-3}%
{\footnotesize 
Invariant:} & 
{\footnotesize 
Homology group,} $H_{C}\left(  M\right)  $ & 
{\footnotesize 
Fukaya category,}
$Fuk\left(  \mathcal{M}^{SLag}\left(  X\right)  \right)  .$ & $%
\begin{array}
[c]{l}%
\,\\
\,
\end{array}
$\\\cline{1-3}%
\end{tabular}
\]

These associations can be formalized to form a TQFT \cite{At TQFT}. Namely we
associate an additive category $F\left(  X\right)  =Fuk\left(  \mathcal{M}%
^{SLag}\left(  X\right)  \right)  $ to a closed almost Calabi-Yau threefold
$X$, a functor $F\left(  M\right)  :F\left(  X_{0}\right)  \rightarrow
F\left(  X_{1}\right)  $ to an almost $G_{2}$-manifold $M$ with asymptotically
cylindrical ends $X_{1}-X_{0}=X_{1}\cup\bar{X}_{0}$. They satisfy
\[%
\begin{tabular}
[c]{ll}%
(i) & $F\left(  \phi\right)  =\text{the additive tensor category of vector
spaces }\left(  \left(  Vec\right)  \right)  \text{,}$\\
(ii) & $F\left(  X_{1}\amalg X_{2}\right)  =F\left(  X_{1}\right)  \otimes
F\left(  X_{2}\right)  .$%
\end{tabular}
\]
For example, when $M$ is a closed $G_{2}$-manifold, that is a cobordism
between empty manifolds, then we have $F\left(  M\right)  :\left(  \left(
Vec\right)  \right)  \rightarrow\left(  \left(  Vec\right)  \right)  $ and the
image of the trivial bundle is our homology group $H_{C}\left(  M\right)  $.

\section{More TQFTs}

Notice that all TQFTs we propose in this paper are formal mathematical
constructions. Besides the lack of compactness for the moduli spaces, the
\textit{obstruction} issue is also a big problem if we try to make these
theories rigorous. This problem is explained to the author by a referee. 

There are other TQFTs naturally associated to Calabi-Yau threefolds and
$G_{2}$-manifolds but (1) they do not involve nontrivial coupling between
submanifolds and bundles and (2) new difficulties arise because of
corresponding moduli spaces for Calabi-Yau threefolds have virtual dimension
zero and could be singular. They are essentially in the paper by Donaldson and
Thomas \cite{DT}.

\bigskip

\textbf{TQFT of associative cycles}

We assume that $M$ is a $G_{2}$-manifold, i.e. $\Omega$ is parallel rather
than closed. Three dimensional submanifolds $A$ in $M$ calibrated by $\Omega$
is called \textit{associative submanifolds} and they can be characterized by
$\chi|_{A}=0$ (\cite{HL}) where $\chi\in\Omega^{3}\left(  M,T_{M}\right)  $ is
defined by $\left\langle w,\chi\left(  x,y,z\right)  \right\rangle =\ast
\Omega\left(  w,x,y,z\right)  $. We define a \textit{parametrized A-cycle} to
be a pair $\left(  f,D_{E}\right)  \in\mathcal{C}_{A}=Map\left(  A,M\right)
\times\mathcal{A}\left(  E\right)  ,$ with $f:A\rightarrow M$ a parametrized
A-submanifold and $D_{E}$ is a unitary flat connection on a Hermitian vector
bundle $E$ over $A$. There is also a natural $\mathcal{G}$-invariant closed
one form $\Phi_{A}$ on $\mathcal{C}_{A}$ given by
\[
\Phi_{A}\left(  f,D_{E}\right)  \left(  v,B\right)  =\int_{A}TrF_{E}\wedge
B+\left\langle f^{\ast}\chi,v\right\rangle _{T_{M}},
\]
for any $\left(  v,B\right)  \in\Gamma\left(  A,f^{\ast}T_{M}\right)
\times\Omega^{1}\left(  A,ad\left(  E\right)  \right)  =T_{\left(
f,D_{E}\right)  }\mathcal{C}_{A}$ $.$ Its zero set is the moduli space of
A-cycles in $M$. As before, we could formally apply arguments in Witten's
Morse theory to $\Phi_{A}$ and define a homology group $H_{A}\left(  M\right)
$.

The corresponding category associated to a Calabi-Yau threefold $X$ would be
the Fukaya-Floer category of the moduli space of unitary flat bundles over
holomorphic curves in $X$, denote $\mathcal{M}^{curve}\left(  X\right)  $. We
summarize these in the following table:
\[%
\begin{tabular}
[c]{|c||c|c|c}\cline{1-3}%
{\footnotesize 
Manifold:} & $G_{2}$
{\footnotesize 
-manifold,} $M^{7}$ & 
{\footnotesize 
CY threefold,} $X^{6}$ & $%
\begin{array}
[c]{l}%
\,\\
\,
\end{array}
$\\\cline{1-3}%
{\footnotesize 
SUSY Cycles:} & 
{\footnotesize 
A-submfds.+ flat bundles} & 
{\footnotesize 
Holomorphic curves+ flat bundles} & $%
\begin{array}
[c]{l}%
\,\\
\,
\end{array}
$\\\cline{1-3}%
{\footnotesize 
Invariant:} & 
{\footnotesize 
Homology group,} $H_{A}\left(  M\right)  $ & 
{\footnotesize 
Fukaya category,}
$Fuk\left(  \mathcal{M}^{curve}\left(  X\right)  \right)  .$ & $%
\begin{array}
[c]{l}%
\,\\
\,
\end{array}
$\\\cline{1-3}%
\end{tabular}
\]

\textbf{TQFT of Donaldson-Thomas bundles}

We assume that $M$ is a seven manifold with a $G_{2}$-structure such that its
positive three form $\Omega$ is co-closed, rather than closed, i.e.
$d\Theta=0$ with $\Theta=\ast\Omega$. In \cite{DT} Donaldson and Thomas
introduce a first order Yang-Mills equation for $G_{2}$-manifolds,
\[
F_{E}\wedge\Theta=0\text{.}%
\]
Their solutions are the zeros of the following gauge invariant one form
$\Phi_{DT}$ on $\mathcal{A}\left(  E\right)  $,%

\[
\Phi_{DT}\left(  D_{E}\right)  \left(  B\right)  =\int_{M}Tr\left[
F_{E}\wedge B\right]  \wedge\Theta\text{,}%
\]
for any $B\in\Omega^{1}\left(  M,ad\left(  E\right)  \right)  =T_{D_{E}%
}\mathcal{A}\left(  E\right)  $. This one form $\Phi_{DT}$ is closed because
of $d\Theta=0$. As before, we can formally define a homology group
$H_{DT}\left(  M\right)  $. The corresponding category associated to a
Calabi-Yau threefold $X$ should be the Fukaya-Floer category of the moduli
space of Hermitian Yang-Mills connections over $X$, denote $\mathcal{M}%
^{HYM}\left(  X\right)  $. Again we summarize these in a table:%

\[%
\begin{tabular}
[c]{|c||c|c|c}\cline{1-3}%
{\footnotesize 
Manifold:} & $G_{2}$
{\footnotesize 
-manifold,} $M^{7}$ & 
{\footnotesize 
CY threefold,} $X^{6}$ & $%
\begin{array}
[c]{l}%
\,\\
\,
\end{array}
$\\\cline{1-3}%
{\footnotesize 
SUSY Cycles:} & 
{\footnotesize 
DT-bundles} & 
{\footnotesize 
Hermitian YM-bundles} & $%
\begin{array}
[c]{l}%
\,\\
\,
\end{array}
$\\\cline{1-3}%
{\footnotesize 
Invariant:} & 
{\footnotesize 
Homology group,} $H_{DT}\left(  M\right)  $ & 
{\footnotesize 
Fukaya category,}
$Fuk\left(  \mathcal{M}^{HYM}\left(  X\right)  \right)  .$ & $%
\begin{array}
[c]{l}%
\,\\
\,
\end{array}
$\\\cline{1-3}%
\end{tabular}
\]

It is an interesting problem to understand the transformations of these TQFTs
under dualities in M-theory.

\bigskip

\textit{Acknowledgments: This paper is partially supported by NSF/DMS-0103355.
The author expresses his gratitude to J.H. Lee, R. Thomas, A. Voronov and X.W.
Wang for useful discussions. The author also thank the referee for many useful comments.}

\bigskip

\end{document}